\documentclass[11pt]{article}
\usepackage{mathrsfs}
\usepackage{}
\usepackage{amsfonts}
\usepackage{amssymb}
\usepackage{graphicx}
\usepackage{amsmath}

\def\sqr#1#2{{\vcenter{\vbox{\hrule height.#2pt
              \hbox{\vrule width.#2pt height#1pt \kern#1pt \vrule
width.#2pt}
              \hrule height.#2pt}}}}
\def\signed #1{{\unskip\nobreak\hfil\penalty50
              \hskip2em\hbox{}\nobreak\hfil#1
              \parfillskip=0pt \finalhyphendemerits=0 \par}}
\def\endpf{\signed {$\sqr69$}}
\def\3n{\negthinspace \negthinspace \negthinspace }
\def\2n{\negthinspace \negthinspace }
\def\1n{\negthinspace }

\def\ms{\medskip}

\def\({\Big (}
\def\){\Big )}
\def\[{\Big[}
\def\]{\Big]}

\def\be{\begin{equation}}
\def\bel{\begin{equation}\label}
\def\ee{\end{equation}}
\def\bea{\begin{eqnarray}}
\def\eea{\end{eqnarray}}
\def\bt{\begin{theorem}}
\def\et{\end{theorem}}
\def\bc{\begin{corollary}}
\def\ec{\end{corollary}}
\def\bl{\begin{lemma}}
\def\el{\end{lemma}}
\def\bp{\begin{proposition}}
\def\ep{\end{proposition}}
\def\br{\begin{remark}}
\def\er{\end{remark}}
\def\ba{\begin{array}}
\def\ea{\end{array}}
\def\bd{\begin{definition}}
\def\ed{\end{definition}}

\newtheorem{lemma}{Lemma}[section]
\newtheorem{remark}{Remark}[section]

\newtheorem{theorem}{Theorem}[section]
\newtheorem{corollary}{Corollary}[section]

\newtheorem{definition}{Definition}[section]
\newtheorem{proposition}{Proposition}[section]

\makeatletter
   
   \@addtoreset{equation}{section}
\makeatother
\allowdisplaybreaks
\begin{document}

\title{\bf A note on damped wave equations with a nonlinear dissipation in non-cylindrical domains \thanks{This work is
supported by the NSF of China under grants 11471070, 11771074 and
11371084.}}

\author{Lingyang Liu\thanks{School of Mathematics and Statistics, Northeast Normal
University, Changchun 130024, China. E-mail address:
liuly938@nenu.edu.cn. \ms }
}

\date{}

\maketitle

\begin{abstract}
In this paper, we study the large time behavior of a class of wave equation with a nonlinear dissipation in non-cylindrical domains. The result we obtained here relaxes the conditions for the nonlinear term coefficients (in precise, that is $\beta(t)|u|^\rho u$) in \cite{alb} and \cite{ha} (which require $\beta(t)$ to be a constant or $\beta(t)$ to be decreasing with time $t$) and has less restriction for the defined regions.
\end{abstract}


\section{Introduction and main results}

Fix $t\geq0.$ Let $\Omega_t$ be a bounded domain in $\mathbb R.$ Given $T>0.$ Set $\widehat{Q}_T=\Omega_t\times(0,T)$ and denote by $\widehat{\Sigma}_T$ the lateral boundary of $\widehat{Q}_T.$ Consider the following wave equation with a nonlinear dissipation in the non-cylindrical domain $\widehat{Q}_T:$
\begin{equation}\label{e1.1}
\left\{\begin{array}{ll}
u^{''}-\Delta u+au'+bu+\beta(t)|u|^\rho u=0& (x,t)\in \widehat{Q}_T,\\[2mm]
u=0 &(x,t)\in\widehat{\Sigma}_T,\\[2mm]
u(x,0)=u_0(x), \ u'(x,0)=u_1(x)&x\in\Omega_0,
\end{array}\right.
\end{equation}
where $(u_0, u_1)$ is any given initial couple, $(u, u')$ is the state variable and $a,b>0.$

In order to study the qualitative theory of (\ref{e1.1}), we need the following assumptions on the domain $\widehat{Q}_T:$

\medskip

\noindent (A1) $\alpha\in C^2[0,T]$ such that $\alpha(0)=1,$ $\alpha'(t)\geq0$ and $\sup\limits_{ t\in[0,T]}\alpha'(t)<1.$

\medskip

\noindent (A2) $\beta(t),\beta'(t)\geq0,$ $t\in[0,T]$ and $\beta'\in L^\infty(0,T).$

\medskip

\noindent (A3) if $n>2,$ then $\displaystyle0<\rho\leq\frac{2}{n-2};$ if $n=1$ or $n=2,$ then $0<\rho<\infty.$

\medskip

The wellposedness result for (\ref{e1.1}) is stated as follows:

\begin{theorem}\label{d1.2}Let $u_0\in H^2_0(0,1)$ and $u_1\in H^1_0(0,1).$ If assumptions (A1)-(A3) hold, then there exists a unique strong solution $u$ of problem (\ref{e1.1}) such that
 $u\in L^\infty\big(0,T;H^1_0(\Omega_t)\cap H^2(\Omega_t)\big),$ $u_t\in L^\infty\big(0,T;H^1(\Omega_t)\big),$ $u_{tt}\in L^\infty\big(0,T;L^2(\Omega_t)\big),$
and
\begin{equation*}
\big(u^{''}-\Delta u+au'+bu+\beta(t)|u|^\rho u, \phi\big)(t)=0, \ a.e.\ t\in(0,T),
\end{equation*}
where $\phi(t)$ is an arbitrary function from $L^2(\mathbb R^1).$ In addition, $u(0)=u_0,$ $u_t(0)=u_1.$
\end{theorem}

The proof of Theorem \ref{d1.2} is quite similar to the proof of wellposedness results in \cite{fer}, so we omit it (but what we need to point out is that since the assumption (A2) is different from $\beta'\leq0,$  the result we obtained here just admits the solution to belong to $L^\infty\big(0,T;H^1_0(\Omega_t)\cap H^2(\Omega_t)\big),$ not to $L^\infty\big(0,\infty;H^1_0(\Omega_t)\cap H^2(\Omega_t)\big)$).

\begin{lemma}[\cite{Ma}]\label{}
Suppose that $\widehat{Q}_T$ has a regular lateral boundary $\widehat{\Sigma}_T.$ If $u\in C^1\big(\mathbb R; L^2(\Omega_t)\big),$ then we have
\begin{eqnarray*}
&&\frac{d}{dt}\int_{\Omega_t}u(x,t)dx=\int_{\Omega_t}\frac{d}{dt}u(x,t)dx+\int_{\Gamma_t}u(x,t)\dot{x}n_xd\sigma\\[2mm]
&&=\int_{\Omega_t}\frac{d}{dt}u(x,t)dx-\int_{\Gamma_t}u(x,t)n_td\sigma,
\end{eqnarray*}
where $\Gamma_t$ is the boundary of $\Omega_t,$  $\dot{x}$ is the velocity of $x\in\Gamma_t,$ and  $n=(n_x,n_t)$ is the unit exterior normal to $\widehat{\Sigma}_T.$
Moreover, it was observed that for $u\in H^1(\widehat{Q}_T)$ with $u=0$ on $\widehat{\Sigma}_T$ (all tangential derivative of $u$ also vanishes on $\widehat{\Sigma}_T$), Consequently the full gradient of $u$ satisfies $\nabla_{x,t}u=(\partial_n u)n$ which implies that
\begin{equation*}
  u_t=(\partial_n u) n_t      \quad \mbox{and} \quad \nabla_xu=(\partial_n u) n_x.
\end{equation*}
\end{lemma}

The energy of system (\ref{e1.1}) $\mathscr{E}(t)$ is given by
\begin{equation*}
 \mathscr{E}(t)= \int_{\Omega_t}\Big[\frac{1}{2}u_t^2(t)+\frac{1}{2}u_x^2(t)+\frac{1}{2}u^2(t)+\beta(t)\frac{1}{\rho+2}|u(t)|^{\rho+2}\Big]dx.
\end{equation*}

Then the main result of this paper is stated as follows.
\begin{theorem}\label{}
One can find $\lambda>0$ and $\beta(t)$ satisfying $\lambda(\rho+1)\beta(t)\geq\beta'(t),$ such that the inequality
\begin{equation}\label{e1}
  \mathscr{E}(t)\leq C\mathscr{E}(0)\varphi^{-1}(t),
\end{equation}
hold, where $\varphi(t)$ is chosen by $\varphi(t)=e^{\lambda t},$ $C$ is some positive constant.
\end{theorem}

\noindent {\bf Proof.} Firstly, let $\varphi$ be a unknown continuous function. Secondly, Multiplying both sides of the first equation in (\ref{e1.1}) by $(u_t+\lambda u)\varphi(t),$ where $\lambda>0,$ and then integrating it on $(0,T)\times \Omega_t,$ we get
\begin{equation*}
 \displaystyle \int^T_0\int_{\Omega_t} \big(u^{''}-\Delta u+au'+bu+\beta(t)|u|^\rho u\big)(u_t+\lambda u)\varphi(t)dxdt=0.
\end{equation*}
Calculating the above equality, we have
\begin{eqnarray}\label{e2}
\begin{array}{rl}
   &\displaystyle\int^T_0\int_{\Omega_t}u^{''}(u_t+\lambda u)\varphi(t)dxdt\\[6mm]
  =&\!\!\!\displaystyle\int^T_0\int_{\Omega_t}\Big[\big(\frac{1}{2}u_t^2\varphi(t)\big)_t+\big(\lambda\varphi(t)uu_t\big)_t-\lambda\varphi(t)u_t^2-\lambda\varphi'(t)uu_t-\frac{1}{2}\varphi'(t)u_t^2\Big]dxdt\\[6mm]
  =&\!\!\!\displaystyle\int_{\Omega_T}\big(\frac{1}{2}u_t^2(T)\varphi(T)+\lambda\varphi(T)u(T)u_t(T)\big)dx-\int_{\Omega_0}\big(\frac{1}{2}u_t^2(0)\varphi(0)+\lambda\varphi(0)u(0)u_t(0)\big)dx\\[6mm]
&\!\!\!+\displaystyle\int^T_0\int_{\Gamma_t}\frac{1}{2}u_t^2\varphi(t)n_td\sigma dt-\int^T_0\int_{\Omega_t}\big[\lambda\varphi(t)u_t^2+\lambda\varphi'(t)uu_t+\frac{1}{2}\varphi'(t)u_t^2\big]dxdt,
\end{array}
\end{eqnarray}

\medskip

\begin{eqnarray}\label{e3}
\begin{array}{rl}
&\displaystyle\int^T_0\int_{\Omega_t}-\Delta u(u_t+\lambda u)\varphi(t)dxdt\\[6mm]
=&\!\!\!\displaystyle\int^T_0\int_{\Omega_t}\Big[\big(-u_xu_t\varphi(t)\big)_x+u_xu_{tx}\varphi(t)-\big(u_x\lambda u\varphi(t)\big)_x-\lambda\varphi(t)u_x^2dxdt\Big]\\[6mm]
=&\!\!\!\displaystyle\int^T_0\int_{\Omega_t}\Big[\big(-u_xu_t\varphi(t)\big)_x+\big(\frac{1}{2}u_x^2\varphi(t)\big)_t-\frac{1}{2}\varphi'(t)u_x^2-\big(\lambda\varphi(t)uu_x\big)_x+\lambda\varphi(t)u_x^2\Big]dxdt\\[6mm]
=&\!\!\!\displaystyle\int^T_0\int_{\Omega_t}\big(-u_xu_t\varphi(t)\big)_xdxdt+\int_{\Omega_T}\frac{1}{2}u_x^2(T)\varphi(T)dx-\int_{\Omega_0}\frac{1}{2}u_x^2(0)\varphi(0)dx\\[6mm]
&\!\!\!\displaystyle+\int^T_0\int_{\Gamma_t}\frac{1}{2}u_x^2\varphi(t)n_td\sigma dt-\int^T_0\int_{\Omega_t}\big[\frac{1}{2}\varphi'(t)u_x^2-\lambda\varphi(t)u_x^2\big]dxdt,
\end{array}
\end{eqnarray}

\medskip

\begin{equation}\label{e4}
 \int^T_0\int_{\Omega_t} au'(u_t+\lambda u)\varphi(t)dxdt=\int^T_0\int_{\Omega_t}\big[a\varphi(t)u_t^2+a\lambda uu_t\varphi(t)\big]dxdt,
\end{equation}

\medskip

\begin{eqnarray}\label{e5}
\begin{array}{rl}
&\displaystyle\int^T_0\int_{\Omega_t}bu(u_t+\lambda u)\varphi(t)dxdt\\[6mm]
=&\!\!\!\displaystyle\int^T_0\int_{\Omega_t}\big[buu_t\varphi(t)+b\lambda\varphi(t)u^2\big]dxdt\\[6mm]
=&\!\!\!\displaystyle\int^T_0\int_{\Omega_t}\Big[\big(\frac{1}{2}bu^2\varphi(t)\big)_t-\frac{b}{2}\varphi'(t)u^2+b\lambda\varphi(t)u^2\Big]dxdt\\[6mm]
=&\!\!\!\displaystyle\int_{\Omega_T}\frac{1}{2}b\varphi(T)u^2(T)dx-\int_{\Omega_0}\frac{1}{2}b\varphi(0)u^2(0)dx-\int^T_0\int_{\Omega_t}\big[\frac{b}{2}\varphi'(t)u^2-b\lambda\varphi(t)u^2\big]dxdt,
\end{array}
\end{eqnarray}

\medskip

\begin{eqnarray}\label{e6}
&&\!\!\!\!\!\!\!\!\int^T_0\int_{\Omega_t}\beta(t)|u|^\rho u(u_t+\lambda u)\varphi(t)dxdt\notag\\[2mm]
=&&\!\!\!\!\!\!\!\!\int^T_0\int_{\Omega_t}\Big[\beta(t)\big(\frac{1}{\rho+2}|u|^{\rho+2}\big)_t\varphi(t)+\lambda\beta(t)|u|^{\rho+2}\varphi(t)\Big]dxdt\notag\\[2mm]
=&&\!\!\!\!\!\!\!\!\int^T_0\int_{\Omega_t}\Big(\frac{1}{\rho+2}|u|^{\rho+2}\beta(t)\varphi(t)\big)_t-\beta'(t)\varphi(t)\frac{1}{\rho+2}|u|^{\rho+2}-\beta(t)\varphi'(t)\frac{1}{\rho+2}|u|^{\rho+2}\Big]dxdt\notag\\[2mm]
&&\!\!\!\!\!\!\!\!+\int^T_0\int_{\Omega_t}\lambda\beta(t)|u|^{\rho+2}\varphi(t)dxdt\\[2mm]
=&&\!\!\!\!\!\!\!\!\displaystyle\int_{\Omega_T}\beta(T)\varphi(T)\frac{1}{\rho+2}|u(T)|^{\rho+2}dx-\int_{\Omega_0}\beta(0)\varphi(0)\frac{1}{\rho+2}|u(0)|^{\rho+2}dx\notag\\[2mm]
&&\!\!\!\!\!\!\!\!\displaystyle+\int^T_0\int_{\Omega_t}\big[\beta'(t)\varphi(t)\frac{1}{\rho+2}|u|^{\rho+2}+\beta(t)\varphi'(t)\frac{1}{\rho+2}|u|^{\rho+2}-\lambda\beta(t)|u|^{\rho+2}\varphi(t)\big]dxdt.\notag
\end{eqnarray}

\medskip

Adding (\ref{e2}) to (\ref{e6}), we obtain
\begin{eqnarray}\label{plu}
0=\!\!\!\!\!\!\!\!&&\displaystyle\int_{\Omega_T}\big(\frac{1}{2}u_t^2(T)\varphi(T)+\lambda\varphi(T)u(T)u_t(T)\big)dx-\int_{\Omega_0}\big(\frac{1}{2}u_t^2(0)\varphi(0)+\lambda\varphi(0)u(0)u_t(0)\big)dx\notag\\[3mm]
&&\displaystyle+\int^T_0\int_{\Gamma_t}\frac{1}{2}u_t^2\varphi(t)n_td\sigma dt-\int^T_0\int_{\Omega_t}\big[\lambda\varphi(t)u_t^2+\lambda\varphi'(t)uu_t+\frac{1}{2}\varphi'(t)u_t^2\big]dxdt\notag\\[3mm]
&&\displaystyle-\int^T_0\int_{\Omega_t}\big(u_xu_t\varphi(t)\big)_x dxdt+\int_{\Omega_T}\frac{1}{2}u_x^2(T)\varphi(T)dx-\int_{\Omega_0}\frac{1}{2}u_x^2(0)\varphi(0)dx\notag\\[3mm]
&&\displaystyle+\int^T_0\int_{\Gamma_t}\frac{1}{2}u_x^2\varphi(t)n_td\sigma dt-\int^T_0\int_{\Omega_t}\big[\frac{1}{2}\varphi'(t)u_x^2-\lambda\varphi(t)u_x^2\big]dxdt\notag\\[3mm]
&&\displaystyle +\int^T_0\int_{\Omega_t}\big[a\varphi(t)u_t^2+a\lambda uu_t\varphi(t)\big]dxdt\\[3mm]
&&\displaystyle+\int_{\Omega_T}\frac{1}{2}b\varphi(T)u^2(T)dx-\int_{\Omega_0}\frac{1}{2}b\varphi(0)u^2(0)dx-\int^T_0\int_{\Omega_t}\big[\frac{b}{2}\varphi'(t)u^2-b\lambda\varphi(t)u^2\big]dxdt\notag\\[3mm]
&&\displaystyle+\int_{\Omega_T}\beta(T)\varphi(T)\frac{1}{\rho+2}|u(T)|^{\rho+2}dx-\int_{\Omega_0}\beta(0)\varphi(0)\frac{1}{\rho+2}|u(0)|^{\rho+2}dx\notag\\[3mm]
&&\displaystyle+\int^T_0\int_{\Omega_t}\big[-\beta'(t)\varphi(t)\frac{1}{\rho+2}|u|^{\rho+2}-\beta(t)\varphi'(t)\frac{1}{\rho+2}|u|^{\rho+2}+\lambda\beta(t)|u|^{\rho+2}\varphi(t)\big]dxdt.\notag
\end{eqnarray}
Since the assumption (A1) means that

\medskip

\noindent (H1) The domain $\widehat{Q}_T$ is time-like, i.e., $|n_t|<|n_x|.$

\medskip

\noindent (H2) $\widehat{Q}_T$ is monotone increasing, i.e., $\Omega_t$ is expanding with respect to $t$ or $n_t\leq0.$

\begin{eqnarray*}
&&\!\!\!\!\!\!\!\!\int^T_0\int_{\Gamma_t}\big[\frac{1}{2}u_t^2\varphi(t)n_t+\frac{1}{2}u_x^2\varphi(t)n_t\big]d\sigma dt-\int^T_0\int_{\Omega_t}\big(u_xu_t\varphi(t)\big)_xdxdt\\[3mm]
=&&\!\!\!\!\!\!\!\!\int^T_0\int_{\Gamma_t}\big[\frac{1}{2}u_t^2\varphi(t)n_t+\frac{1}{2}u_x^2\varphi(t)n_t\big]d\sigma dt-\int^T_0\int_{\Gamma_t}u_xu_t\varphi(t)n_xd\sigma dt\\[3mm]
=&&\!\!\!\!\!\!\!\!\int^T_0\int_{\Gamma_t} \frac{1}{2}\varphi(t)|\partial_n u|^2(n_t^2-n_x^2) n_td\sigma dt\geq0.
\end{eqnarray*}
Furthermore, (\ref{plu}) yields
\begin{eqnarray}\label{ine}
&&\!\!\!\!\!\!\!\!\displaystyle\int_{\Omega_T}\Big[\frac{1}{2}u_t^2(T)+\lambda u(T)u_t(T)+\frac{1}{2}u_x^2(T)+\frac{1}{2}bu^2(T)+\beta(T)\frac{1}{\rho+2}|u(T)|^{\rho+2}\Big]\varphi(T)dx\notag\\[3mm]
\leq&&\!\!\!\!\!\!\!\!\displaystyle\int_{\Omega_0}\Big[\frac{1}{2}u_t^2(0)+\lambda u(0)u_t(0)+\frac{1}{2}u_x^2(0)+\frac{1}{2}bu^2(0)+\beta(0)\frac{1}{\rho+2}|u(0)|^{\rho+2}\Big]\varphi(0)dx\notag\\[3mm]
&&\!\!\!\!\!\!\!\!\displaystyle+\int^T_0\int_{\Omega_t}\big[\lambda\varphi(t)u_t^2+\lambda\varphi'(t)uu_t+\frac{1}{2}\varphi'(t)u_t^2\big]dxdt+\int^T_0\int_{\Omega_t}\big[\frac{1}{2}\varphi'(t)u_x^2-\lambda\varphi(t)u_x^2\big]dxdt\notag\\[3mm]
&&\!\!\!\!\!\!\!\!\displaystyle-\int^T_0\int_{\Omega_t}\big[a\varphi(t)u_t^2+a\lambda uu_t\varphi(t)\big]dxdt+\int^T_0\int_{\Omega_t}\big[\frac{b}{2}\varphi'(t)u^2-b\lambda\varphi(t)u^2\big]dxdt\notag\\[3mm]
&&\!\!\!\!\!\!\!\!\displaystyle+\int^T_0\int_{\Omega_t}\big[\beta'(t)\varphi(t)\frac{1}{\rho+2}|u|^{\rho+2}+\beta(t)\varphi'(t)\frac{1}{\rho+2}|u|^{\rho+2}-\lambda\beta(t)|u|^{\rho+2}\varphi(t)\big]dxdt.
\end{eqnarray}

We can choose $\varphi(t)=e^{st},$ $s>0.$ In particular, let $\varphi(t)=e^{\lambda t}$ ($\lambda$ be small) and
\begin{eqnarray}\label{bet}
\lambda(\rho+1)\beta(t)\geq\beta'(t).
\end{eqnarray}
We can put
\begin{equation*}
  \beta(t)=e^{\mu t} \quad \mbox{with} \quad  \mu\leq\lambda(\rho+1),
\end{equation*}
or
\begin{equation*}
  \beta(t)=a_nt^n+a_{n-1}t^{n-1}+\cdots+a_1t+a_0,
\end{equation*}
with $a_i>0 (i=0,\cdots,n)$ such that (\ref{bet}) holds.

\medskip

Then the last three terms of inequality (\ref{ine}) are negative. Hence, we deduce
\begin{eqnarray*}
&&\!\!\!\!\!\!\!\!\displaystyle\int_{\Omega_T}\Big[\frac{1}{2}u_t^2(T)+\lambda u(T)u_t(T)+\frac{1}{2}u_x^2(T)+\frac{1}{2}bu^2(T)+\beta(T)\frac{1}{\rho+2}|u(T)|^{\rho+2}\Big]\varphi(T)dx\\[3mm]
\leq&&\!\!\!\!\!\!\!\!\displaystyle\int_{\Omega_0}\Big[\frac{1}{2}u_t^2(0)+\lambda u(0)u_t(0)+\frac{1}{2}u_x^2(0)+\frac{1}{2}bu^2(0)+\beta(0)\frac{1}{\rho+2}|u(0)|^{\rho+2}\Big]\varphi(0)dx.
\end{eqnarray*}
From the above inequality, we finally derive
\begin{equation*}\label{con}
  \mathscr{E}(t)\leq C\mathscr{E}(0)\varphi^{-1}(t),
\end{equation*}
for some constant $C>0.$

\endpf

\medskip

\begin{remark}
If $b=0$ in (\ref{e1.1}), then use the method before, (\ref{ine}) becomes
\begin{eqnarray*}
&&\!\!\!\!\!\!\!\!\displaystyle\int_{\Omega_T}\Big[\frac{1}{2}u_t^2(T)+\lambda u(T)u_t(T)+\frac{1}{2}u_x^2(T)+\beta(T)\frac{1}{\rho+2}|u(T)|^{\rho+2}\Big]\varphi(T)dx\\[3mm]
\leq&&\!\!\!\!\!\!\!\!\displaystyle\int_{\Omega_0}\Big[\frac{1}{2}u_t^2(0)+\lambda u(0)u_t(0)+\frac{1}{2}u_x^2(0)+\beta(0)\frac{1}{\rho+2}|u(0)|^{\rho+2}\Big]\varphi(0)dx\\[3mm]
&&\!\!\!\!\!\!\!\!\displaystyle+\int^T_0\int_{\Omega_t}\big[\lambda\varphi(t)u_t^2+\lambda\varphi'(t)uu_t+\frac{1}{2}\varphi'(t)u_t^2\big]dxdt+\int^T_0\int_{\Omega_t}\big[\frac{1}{2}\varphi'(t)u_x^2-\lambda\varphi(t)u_x^2\big]dxdt\\[3mm]
&&\!\!\!\!\!\!\!\!\displaystyle-\int^T_0\int_{\Omega_t}\big[a\varphi(t)u_t^2+a\lambda uu_t\varphi(t)\big]dxdt\\[3mm]
&&\!\!\!\!\!\!\!\!\displaystyle+\int^T_0\int_{\Omega_t}\big[\beta'(t)\varphi(t)\frac{1}{\rho+2}|u|^{\rho+2}+\beta(t)\varphi'(t)\frac{1}{\rho+2}|u|^{\rho+2}-\lambda\beta(t)|u|^{\rho+2}\varphi(t)\big]dxdt.\\[3mm]
\end{eqnarray*}
In this case, in order to absorb the mixed term $\int^T_0\int_{\Omega_t}a\lambda uu_t\varphi(t)dxdt,$ we must use poincar\'{e} inequality whose coefficients depend on geometry of the domain. That is
\begin{equation*}
\displaystyle \int_{\Omega_t} u^2(x,t)dx\leq |\Omega_t|^2\int_{\Omega_t} u_x^2(x,t)dx.
\end{equation*}
Thus
\begin{eqnarray*}
&&\displaystyle\int^T_0\int_{\Omega_t}a\lambda uu_t\varphi(t)dxdt\leq\int^T_0\int_{\Omega_t}\frac{1}{2}a\lambda^2 \varphi(t)u^2dxdt+\int^T_0\int_{\Omega_t}\frac{1}{2}a\varphi(t)u_t^2dxdt\\[3mm]
&&\leq\int^T_0\int_{\Omega_t}\frac{1}{2}a\lambda^2 |\Omega_t|^2\varphi(t)u_x^2dxdt+\int^T_0\int_{\Omega_t}\frac{1}{2}a\varphi(t)u_t^2dxdt.
\end{eqnarray*}

When $\alpha\in L^\infty(0,\infty),$ and there exist two bounded domains $\Omega_*, \Omega^* \subset\mathbb R^1$ such that $\Omega_*\subset\Omega_\tau\subset\Omega_t\subset\Omega^*, \forall\tau<t.$ Then we have $|\Omega_t|\leq|\Omega^*|,\forall t>0.$ Let $\displaystyle a\lambda|\Omega^*|^2<1.$ With a similar argument as before, we get
\begin{equation*}
  \mathscr{E}(t)\leq C\mathscr{E}(0)\varphi^{-1}(t),\quad t>0,
\end{equation*}
for some constant $C>0.$

If non-cylindrical domains become unbounded in some $X_1$-direction of space, as the time $t$ goes to infinite, and are bounded in other $X_2$-direction of space. Since the projection of it in $X_2$-direction is a bounded open set, written as $w,$ then the Poincar\'{e} inequality in $X_2$-direction turns out
\begin{equation*}
  \displaystyle \int_{\Omega_t} u^2(x,t)dx\leq C_{w}^2\int_{\Omega_t} |\nabla_{X_2}u(x,t)|^2dx\leq C_{w}^2\int_{\Omega_t}|\nabla u(x,t)|^2dx,\quad
\end{equation*}
where $C_{w}$ is the Poincar\'{e} constant.

\medskip

Therefore, the above conclusion is still valid for this case.
\end{remark}

\begin{remark}
For the case of domains becoming unbounded in every spatial direction, as the time $t$ goes to infinite, the condition $b\neq0$ is needed to make (\ref{e1}) true. Otherwise, for any given $T>0,$ let $\lambda=\lambda(T)$ (depending on time $T$) be small and then it follows that
\begin{equation*}
 \mathscr{E}(t)\leq C\mathscr{E}(0)\varphi_T^{-1}(t),\quad 0<t<T,
\end{equation*}
where $\varphi_T^{-1}(t)=e^{-\lambda(T)t}.$

Since Poincar\'{e} inequality does not hold for a fixed number in any totally unbounded area, it seems difficult for us to get an estimate (\ref{e1}) without compensation ($b=0$) and this is also an open problem that has been mentioned in some literature such as \cite{ha}.
\end{remark}

\end{document}